\documentclass[a4paper,10pt]{article}
\usepackage[centertags]{amsmath}
\usepackage{amsfonts}
\usepackage{amssymb}
\usepackage{amsthm}
\usepackage{epsfig}
\usepackage{setspace}
\usepackage{ae}
\usepackage{eucal}
\usepackage[usenames]{color}%
\usepackage{graphicx}

\theoremstyle{plain}
\newtheorem{thm}{Theorem}[section]

\newtheorem{cor}[thm]{Corollary}

\theoremstyle{definition}

\theoremstyle{remark}

\newtheorem{rem}[thm]{Remark}

\title{Singular solutions of Navier Stokes equations with time-dependent external force terms in $L^2$}
\author{J\"org Kampen }

\begin{document}

\maketitle

\begin{abstract}
It is shown  that Navier Stokes equation models of dimension $n\geq 3$ with time dependent external forces in $L^2$ can have singular solutions. 
\end{abstract}

\section{Introduction}

We consider incompressible Navier Stokes equation models of the form 
\begin{equation}\label{Navleray}
\left\lbrace \begin{array}{ll}
\frac{\partial v_i}{\partial t}-\nu\sum_{j=1}^n \frac{\partial^2 v_i}{\partial x_j^2} 
+\sum_{j=1}^n v_j\frac{\partial v_i}{\partial x_j}=f_i\\
\\ \hspace{1cm}+\int_{{\mathbb R}^n}\left( \frac{\partial}{\partial x_i}K_n(x-y)\right) \sum_{j,m=1}^n\left( \frac{\partial v_m}{\partial x_j}\frac{\partial v_j}{\partial x_m}\right) (t,y)dy,\\
\\
\mathbf{v}(0,.)=\mathbf{h},
\end{array}\right.
\end{equation}
to be solved for the velocity $\mathbf{v}=\left(v_1,\cdots ,v_n \right)^T$ on the domain $\left[0,\infty \right)\times {\mathbb R}^n$, where the symbol ${\mathbb R}$ denotes the field of real numbers. The following construction of singular solutions works for dimension $n \geq 3$ due to some structural properties of the equation such as  the less restrictive role of incompressibility if compared to the vorticity equation in dimension less than three. The considerations can also be applied (with slight modifications) if the problem is posed on the domain $\left[0,\infty \right)\times {\mathbb T}^n$, where ${\mathbb T}^n$ is a torus of dimension $n\geq 3$ (equivalent to periodic boundary conditions).  Note that the equation in (\ref{Navleray}) is written in Leray projection form, the constant $\nu >0$ is the viscosity,  $K_n$ is the Laplacian kernel of dimension $n\geq 3$, and $\mathbf{h}=\left(h_1,\cdots ,h_n \right)^T$ are the initial data. Furthermore, the functions $f_i:\left[0,\infty \right)\times {\mathbb R}^n\rightarrow {\mathbb R},~1\leq i\leq n$ are force terms, which are in $L^2$ (at least).


In case $n\geq 3$ we show that such models can have singular solutions for some time-dependent weak external forces $f_i,~1\leq i\leq n$ which are just in $L^2$.

\section{Singular solution construction for $L^2$-data of external forces}
In \cite{web} there has been a discussion of a proposed strong solution of (essentially) a Navier Stokes equation model on the $3$-dimensional torus  with force terms which are just in $L^2$. We show that there are singular solutions for simple types of these models, implying that a global regular existence and uniqueness result cannot be obtained for this general class of models. As remarked, we consider the whole domain with spatial part ${\mathbb R}^n$ here, where an analogous reasoning is possible if the spatial part of the domain is a three-dimensional torus.   
We observe that time dependent external forces of low regularity, e.g., $f_i\in L^2$ for all $1\leq i\leq n$ can consume the damping by the viscosity term such that we have some analogy with the incompressible Euler equation. We denote velocity component functions of the incompressible Euler equation by $v^E_i,~1\leq i\leq n$. The domain of regular existence of the constructed singular solution functions will be $[0,\rho)\times {\mathbb R}^n$, where $\rho>0$ is a positive real number. We construct singularities at the tip of a cone, where some features of the Euler equation can be transferred to the incompressible Navier Stokes equation with weak force terms. Next we implement this idea. For the convenience of the reader we derive the transformation rather explicitly.
First we define functions $w^E_i,~1\leq i\leq n$ in terms of the velocity components $v^E_i,~1\leq i\leq n$ (of a solution of the incompressible Euler equation) by
\begin{equation}\label{veuler}
\frac{w^E_i(s,y)}{\rho-t}=v^E_i(t,x),~s=\frac{t}{\sqrt{\rho^2-t^2}},~y_i=(\rho-t)\arctan(x_i),~1\leq i\leq n.
 \end{equation}
Note that 
\begin{equation}
 t=t(s)=\frac{\rho s}{\sqrt{1+s^2}},~~\mbox{and}~~\frac{ds}{dt}=\frac{\rho^2}{\sqrt{\rho^2-t^2}^3}.
\end{equation}
We also use the abbreviations 
\begin{equation}\label{yiabb}
y_i=y_i(t,x_i)=(\rho-t)\arctan(x_i),~ 1\leq i\leq n
\end{equation} 
in the following if we want to emphasize the dependence of the $y_i$-coordinates on original time- and spatial coordinates. For the initial data we note that
\begin{equation}
w^E_i(0,.)=h^{\rho}_i(.),
\end{equation}
where $h^{\rho}_i(.),~1\leq i\leq n$ denote the transformed data, i.e., $h^{\rho}_i(y)=\rho h_i(x)=v^E_i(0,x)$ for all $1\leq i\leq n$.
The idea of the following singular solution construction is that for some small $\rho >0$ and certain regular data $h^{\rho}_i,~1\leq i\leq n$ corresponding to data $h_i,~1\leq i\leq n$ in original coordinates with strong polynomial decay at spatial infinity  we have global solution branches in time $s$-coordinates
 \begin{equation}
 w^E_i:K_{\rho}\rightarrow {\mathbb R},~1\leq i\leq n,~\mbox{where }~K_{\rho}~\mbox{is a cone, cf. (\ref{cone}), }
 \end{equation}
which correspond to local-time solutions on a small time interval $[0,\rho)$ in original time $t$-coordinates, where for some regular data $h^{\rho}_i,~1\leq i\leq n$ we have in addition
\begin{equation}
w^E_i(0,0)=h^{\rho}_i(0)\not = 0,\lim_{t(s)\uparrow \rho}w^E_i(s,0)\neq 0,~\mbox{for some }i_0.
\end{equation}
Here the domain $K_{\rho}$ is the image of the domain $[0,\rho)\times {\mathbb R}^n$ under the transformation $(t,x)\rightarrow (s,y)$ in (\ref{veuler}). The solutions are viscosity limits of local solutions of a family of related equations closely related to incompressible Navier Stokes equations (a family parameterized by the viscosity $\nu >0$). The function $w^E_i, 1\leq i\leq n$ is well defined at the tip of an infinite cone corresponding to the point $(\rho,0)$ in original coordinates, and in addition we can prove that for some $\rho >0$ and for $s=s(t)\uparrow \infty$ or $t\uparrow \rho$ we have that $\lim_{t\uparrow \rho }w_i(s(t),0)\neq 0$ exists  for some $1\leq i\leq n$. This corresponds then to a singularity at the tip of a cone $(\rho,0)$ of the corresponding velocity function component $v^E_i$ evaluated at that point, since
\begin{equation}
v^E_i(t,0)=\frac{w^E_i(s(t),0)}{\rho-t}~\mbox{ for all }~t<\rho.
\end{equation}
This construction is only possible as the transformed function $w^E_i,~1\leq i\leq n$ is supported on a spatially finite cone (the support of the cone is finite with respect to the spatial variables as we use the $\arctan$-function in the transformation). 
Note that for the change of the spatial measure we have
\begin{equation}\label{spatmeas}
dy_1dy_2\cdots dy_n=(\rho-t)^n\Pi_{i=1}^n\frac{dx_i}{1+x_i^2}.
\end{equation}
We shall observe that a local contraction result for the function $w^E_i,~1\leq i\leq n$ can be obtained via an iterative fixed point scheme of viscosity approximations $w^{\nu,E}_i,~1\leq i\leq n$, where the latter function essentially satisfies a related equation with additional viscosity term $-\nu\Delta w^{\nu,E}$ on the left side of the equation ($\nu>0$ is a diffusion constant).  
The measure in (\ref{spatmeas}) ensures that the potential damping term created by a transformation as in (\ref{veuler}) of the form
\begin{equation}\label{damp}
-\frac{\sqrt{\rho^2-t^2}^3}{\rho^2(\rho-t)}w^E_i,
\end{equation}
can be integrated in $(s,y)$-time coordinates globally. Note that otherwise, the condition $\lim_{t\uparrow \rho }w_i(s(t),0)\neq 0$ can not be guaranteed. 
Moreover, the latter term turns out to be relatively strong compared to the nonlinear Euler terms for this transformation on the small time interval $[0,\rho]$ (cf. below). Note that without a spatial change of the measure the coefficient $\frac{\sqrt{\rho^2-t^2}^3}{\rho^2(\rho-t)}$ of the potential term integrates over $t\in [0,\rho]$ corresponding to $s\in \left[ 0,\frac{\rho-\epsilon}{\sqrt{\rho^2-(\rho-\epsilon)^2}}\right] $ for $0<\epsilon <\rho$ as (cf. (\ref{st}) below)
\begin{equation}
\begin{array}{ll}
\int_{t(s)\in [0,\rho-\epsilon]}\frac{\sqrt{\rho^2-t(s)^2}^3}{\rho^2(\rho-t(s))}ds
=\int_0^{\rho-\epsilon}\frac{\sqrt{\rho^2-t^2}^3}{\rho^2(\rho-t)}\frac{ds}{dt}dt\\
\\
=\int_0^{\rho-\epsilon}\frac{1}{(\rho-t)}dt=-\ln(\rho-(\rho-\epsilon))+\ln(\rho)=-\ln\left(\epsilon\right)+\ln(\rho) \uparrow \infty
\end{array}
\end{equation} 
as $\epsilon$ goes to $0$ and $\rho >0$ is fixed. Hence we need a spatial transformation which ensures that the volume of the cone of support is small enough such that $w^E_i(s(t),0)$ becomes well defined as $t\uparrow \rho$ while it is different from zero for at least one index $1\leq i\leq n$.

Next we derive the  transformation in detail for convenience of the reader. Note that for given $t$ the transformed spatial variable $y_i\in \left]-(\rho-t)\frac{\pi}{2},(\rho-t)\frac{\pi}{2}\right[$ corresponds to $x_i\in {\mathbb R}$. More precisely, the transformed functions are supported on a cone (considering the closure of the spatial interval at each time) 
\begin{equation}\label{cone}
K_{\rho}:=\left\lbrace (s,y)|0\leq t(s)\leq \rho~\&~-(\rho-t(s))\frac{\pi}{2}\leq y_i\leq (\rho-t(s))\frac{\pi}{2}\right\rbrace. 
\end{equation}
A time section of the cone $K_{\rho}$ at time $s$ is denoted by $K^s_{\rho}$, i.e., we write
\begin{equation}
K^s_{\rho}:=\left\lbrace (\sigma,y)\in K_{\rho}|\sigma=s\right\rbrace. 
\end{equation}
The cone is of infinite 'height' with respect to $s$-coordinates and of height $\rho$ with respect to $t$-coordinates, where the basis (at $s=t=0$) is the cube  $\left]-\rho\frac{\pi}{2},\rho\frac{\pi}{2}\right[^n$ in transformed coordinates.
We get
\begin{equation}
 v^E_{i,t}=
\frac{w^E_i(s,.)}{\left( \rho-t\right)^2}
+\frac{w^E_{i,s}(s,.)}{\rho-t}\frac{ds}{dt}-\sum_{j=1}^nw^E_{,j}(s,y)
\arctan(x_i),
\end{equation}
where
\begin{equation}\label{st}
\begin{array}{ll}
\frac{ds}{dt}=\frac{1}{\sqrt{\rho^2-t^2}}+\frac{-\frac{1}{2}t(-2t)}{ \sqrt{\rho^2-t^2}^3}
=\frac{\rho^2}{\sqrt{\rho^2-t^2}^3}.
\end{array}
\end{equation}
Furthermore,
\begin{equation}\label{firstder}
v^E_{i,j}=\frac{w^E_{i,j}(s,y)}{\rho-t}\frac{d y_j}{d x_j}=\frac{w^E_{i,j}(s,y)}{\rho-t}\frac{\rho-t}{1+x_j^2}=\frac{w^E_{i,j}(s,y)}{1+x_j^2},
\end{equation}
such that Burgers terms transform like $v^E_jv^E_{i,j}=\frac{1}{(\rho-t)(1+x_j^2)}w^E_jw^E_{i,j}$ and Leray data, i.e., data of the Poisson equation for the pressure, transform like $v^E_{j,k}v^E_{k,j}=\frac{1}{(1+x_j^2)(1+x_k^2)}w^E_{j,k}w^E_{k,j}$. Multiplying the inverse of the coefficient of the time derivative (with respect to the time variable $s$) we get an additional factor $\frac{1}{\rho^2}(\rho -t)\sqrt{\rho^2-t^2}^3$ for all these terms.  
Hence, if $v^E_i,~1\leq i\leq n$ satisfies the incompressible Euler equation on the interval $[0,\rho)$, then $w^E_i,~1\leq i\leq n$ satisfies
\begin{equation}\label{eulertrans}
\begin{array}{ll}
w^E_{i,s}+\sum_{j=1}^n\frac{\sqrt{\rho^2-t^2}^3}{\rho^2(1+x_j^2)}w^E_j\frac{\partial w^E_i}{\partial y_j}-\sum_{j=1}^n\frac{1}{\rho^2}(\rho -t)\sqrt{\rho^2-t^2}^3\arctan(x_j)w^E_{,j}(s,y)\\
\\
=-\frac{\sqrt{\rho^2-t^2}^3}{\rho^2(\rho-t)}w^E_i
+\int_{K^s_{\rho}}\sum_{j,m=1}^n\frac{(\rho-t)\sqrt{\rho^2-t^2}^3}{\rho^2(1+x_j^2)(1+x_m^2)}\times\\
\\
\times\left( K_{n,i}^*(.-z)\right) \left( \frac{\partial w^E_m}{\partial y_j}\frac{\partial w^E_j}{\partial y_m}\right) (s,z)\frac{\Pi_{i=1}^n(1+x_i^2)}{(\rho-t)^n}dz,
\end{array}
\end{equation}
on the cone $K_{\rho}$, where for $n\geq 3$
\begin{equation}\label{k*}
K_{n,i}(x)=c\frac{x_i}{|x|^n},~\mbox{i.e.}~K_{n,i}^*(y)=c\frac{x_i(s,y_i)}{|(x_1(s,y_1),\cdots,x_n(s,y_n))|^n}
\end{equation}
(with a well-known dimension dependent constant $c$), and where  $x_i=x_i(s,y_i)=\tan\left(\frac{y_i}{\rho-t(s)}\right)$ with $t(s)$ defined as above. Note that the equation in (\ref{eulertrans}) is with respect to $s$ and $y_i$-coordinates, where we understand that $x_i=x_i(s,y_i)$ (and $x_i=x_i(s,z_i)$ in the measure of the Leray projection term) according to the transformation above. The integrals are -strictly speaking- on the domain $K_{\rho}$, but we shall use the convenience of classical representations of solutions of approximating equations which are  supported on the whole domain. 
Note that in the Leray projection term $x_i$ is a function which is convoluted, and the support for the solution is the domain $K_{\rho}$.  We always understand $t=t(s)$ according to the transformation above. Recall that $\frac{\Pi_{i=1}^n(1+x_i^2)}{(\rho-t)^n}dy_1dy_2\cdots dy_n=\Pi_{i=1}^n dx_i.$
The equation in (\ref{eulertrans}) is studied along with a family of incompressible Navier Stokes equations for functions $w^{\nu,E}_i,~1\leq i\leq n$, which solve the same equation as in (\ref{eulertrans}), but with an additional viscosity term
\begin{equation}
-\nu \Delta w^{\nu ,E}_i.
\end{equation}

Note that for any $\rho >0$ the time dependent coefficients are bounded on the interval $[0,\rho]$ which follows from the observation
\begin{equation}\label{coeff1}
\lim_{\rho\downarrow 0}\sup_{t\in [0,\rho]}\frac{\sqrt{\rho^2-t^2}^3}{\rho^2}=0,
\end{equation}
for a coefficient of the Burgers term, the Leray projection term and the artificial convection term, while
\begin{equation}\label{coeff2}
\lim_{\rho\downarrow 0}\sup_{t\in [0,\rho]}\frac{\sqrt{\rho^2-t^2}^3}{\rho^2(\rho-t)}=\lim_{\rho\downarrow 0}\sup_{t\in [0,\rho]}\frac{\sqrt{\rho^2-t^2}(\rho+t)}{\rho^2}=1
\end{equation}
for the potential damping term.
Moreover, simple observations as in (\ref{coeff1}) show that the coefficients become arbitrarily small for small $\rho$, while the observation in (\ref{coeff2}) shows that the coefficient of the damping term is comparatively large (of order $1$ for small $\rho>0$). Although the problem for the function $w^E_i,~1\leq i\leq n$ is defined on a global time interval in $s$-time coordinates we can get global existence via a local contraction argument (for small $\rho$), i.e., a local contraction argument leads to local regular existence in the time interval $t\in \left[0,\rho\right]$ for $\rho>0$ small enough corresponding to the global time interval $s\in [0,\infty)$. In order to ensure local contraction it is useful to have strong polynomial decay at infinity of the data functions $h_i$, i.e., for some $m\geq 2$ we have
\begin{equation}
h_i\in  {\cal C}^{m(n+1)}_{pol,m},~1\leq i\leq n,
\end{equation}
where
\begin{equation}\label{pol1}
{\cal C}^{l}_{pol,m}={\Big \{} f:{\mathbb R}^D\rightarrow {\mathbb R}: \\
\\
\exists c>0~\forall |x|\geq 1~\forall 0\leq |\gamma|\leq m~{\big |}D^{\gamma}_xf(x){\big |}\leq \frac{c}{1+|x|^l} {\Big \}}. 
\end{equation}
\begin{rem} In this article we use the function space ${\cal C}^{m(n+1)}_{pol,m}$ for a local time iterative scheme of approximative solutions for $w^{\nu,E}_i,~1\leq i\leq n$ (and then obtain $w^E_i,~1\leq i\leq n$ in the viscosity limit), where the mass is concentrated almost on a (spatially compact) cone. 
The choice of the strong function space ${\cal C}^{m(n+1)}_{pol,m}$ can illustrate how some iterative schemes 'can be tamed' by consideration of strong function spaces although a convolution with Gaussians may transport mass to high frequencies. This taming holds for local time schemes, where an additional argument shows that it holds even for global time schemes. Consider first local time  iterative solution schemes $v^{(k)}_i,~1\leq i\leq n,~k\geq 0$ of (\ref{Navleray}) with $f_i\equiv 0,~1\leq i\leq n$ in terms of classical representations with the fundamental solution $G_{\nu}$  of the equation $\frac{\partial p}{\partial t}-\nu\sum_{j=1}^n \frac{\partial^2 p}{\partial x_j^2}=0$. We define for all $t\geq 0$ and $x\in {\mathbb R}^n$
\begin{equation}
v^{(0)}_i(t,x):=\int_{{\mathbb R}^n}h_i(y) G_{\nu}(t,x;0,y)dy=h_i\ast_{sp}G_{\nu},
\end{equation}
and for $k\geq 1$
\begin{equation}\label{Navleraysol}
\begin{array}{ll}
v^{(k)}_i
=h_i\ast_{sp}G_{\nu}+\sum_{j=1}^n \left( v^{(k-1)}_j\frac{\partial v^{(k-1)}_i}{\partial x_j}\right) \ast G_{\nu}\\
\\ \hspace{1cm}+\sum_{j,m=1}^n\int_{{\mathbb R}^n}\left( \frac{\partial}{\partial x_i}K_n(.-y)\right) \sum_{j,m=1}^n\left( \frac{\partial v^{(k-1)}_m}{\partial x_j}\frac{\partial v^{(k-1)}_j}{\partial x_m}\right) (.,y)dy\ast G_{\nu},
\end{array}
\end{equation}
where $\ast$ denotes convolution with respect to space and time and $\ast_{sp}$ denotes convolution with respect to the spatial variables. In the computation of the increment
\begin{equation}
\delta v_i:=v_i-h_i\ast_{sp}G_{\nu}
\end{equation}
the effect of the initial data smoothing can be eliminated by the consideration of a related iterative solution scheme. Define
\begin{equation}
\delta v^{(0)}_i(t,x):=0,~ \mbox{and}~\delta v^{(k)}_i:=v^{(k)}_i-v^{(k-1)}_i~ \mbox{for $k\geq 1$}. 
\end{equation}
The nonlinear quadratic terms have the effect that the effect of decreasing the degree of spatial decay at spatial infinity  caused by the convolutions with the Gaussians or first order spatial derivatives of the Gaussian and with first order derivatives of the Laplacian kernel is offset by the effect of quadratic powers of functions with strong spatial decay. You prove straightforwardly that for local time $t\geq 0$ and $m\geq 2$
\begin{equation}
\delta v^{(k)}(t,.)\in {\cal C}^{m(n+1)}_{pol,m},~v^{(k)}(t,.)\in {\cal C}^{m(n+1)-\mu}_{pol,m}
\end{equation}
for some $\mu\in (0,1)$ which accounts for the convolution effects of the linear term (in  this case the convoluted initial data).  For local time you straightforwardly show that the increment functions $\delta v^{(k)}(t,.)$ are uniformly bounded in ${\cal C}^{m(n+1)}_{pol,m}$, and that local time contraction holds. This leads to a local time solution representation 
\begin{equation}\label{vsol}
v_i:=v^{(0)}_i+\sum_{k\geq 1}\delta v^{(k)}_i=h_i\ast_{sp}G_{\nu}+\sum_{k\geq 1}\delta v^{(k)}_i\in {\cal C}^{m(n+1)-\mu}_{pol,m},
\end{equation}
where we have some loss of spatial decay of the solution compared to the initial data (due to the convolution of the data, i.e. due to the first term on the right side of (\ref{vsol}) alone). However, this does not mean that this loss of spatial decay explodes as we extend the scheme based on the semigroup property (cf. next remark).
\end{rem}

\begin{rem}
We mention here that the considerations of the last remark can be extended to global schemes (although the following observation is not needed for the local time argument of this article). We have defined global schemes elsewhere where the spatial dependence is in the $H^m\cap C^m$. However in more general models (such as highly degenerated diffusions which satisfy a H\"{o}rmander condition), the assumption of strong spatial polynomial decay is useful in controlled schemes.  The semigroup property of the Gaussian, the preservation of spatial polynomial decay of the functional increments for local time, combined with an external control, or with local damping of the convoluted initial data term $h_i\ast_{sp}G_{\nu}$ (or $v_i(t_0,.)\ast_{sp}G_{\nu}$ lead to global schemes. We only add some remarks concerning spatial polynomial decay here.
Assume that
\begin{equation}
v_i(t_0,.)=h_i\ast_{sp}G_{\nu}+\sum_{k\geq 1}\delta v^{(k)}_i\in {\cal C}^{m(n+1)-\mu}_{pol,m},
\end{equation} 
has been proved for some $t_0>0$. Furthermore assume that
\begin{equation}\label{assdelta}
\sum_{k\geq 1}\delta v^{(k)}_i(t_0,.)\in {\cal C}^{m(n+1)}_{pol,m}.
\end{equation}
We then define a local time iteration scheme on an interval $[t_0,t_1]$ for $t_1>t_0$. First define
\begin{equation}
v^{t_0,(0)}_i:=v_i(t_0,.)\ast_{sp}G_{\nu}.
\end{equation}
Note that for some time interval $[t_0,t_1]$ and $t\in [t_0,t_1]$ we have for some $\mu\in (0,1)$
\begin{equation}
\begin{array}{ll}
v^{t_0,(0)}_i(t,.)=v_i(t_0,.)\ast_{sp}G_{\nu}(t-t_0,.)\\
\\
=\left( h_i\ast_{sp}G_{\nu}(t_0,.)\right)\ast_{sp} 
G_{\nu}(t-t_0,.)+
{\big (}\sum_{k\geq 1}\delta v^{(k)}_i{\big )}\ast_{sp}G_{\nu}(t-t_0,.)\\
\\
=h_i\ast_{sp}G_{\nu}(t,.)+{\big (}\sum_{k\geq 1}\delta v^{(k)}_i(t_0,.){\big )}\ast_{sp}G_{\nu}(t-t_0)\in {\cal C}^{m(n+1)-\mu}_{pol,m},
\end{array}
\end{equation} 
where for the second summand we may use (\ref{assdelta}).
For $k\geq 1$ define on $[t_0,t_1]$
\begin{equation}\label{Navleraysol}
\begin{array}{ll}
v^{t_0,(k)}_i
=v(t_0,.)\ast_{sp}G_{\nu}+\sum_{j=1}^n \left( v^{t_0,(k-1)}_j\frac{\partial v^{t_0,(k-1)}_i}{\partial x_j}\right) \ast G_{\nu}+\\
\\\sum_{j,m=1}^n\int_{{\mathbb R}^n}\left( K_{n,i}(.-y)\right) \sum_{j,m=1}^n\left( \frac{\partial v^{t_0,(k-1)}_m}{\partial x_j}\frac{\partial v^{t_0,(k-1)}_j}{\partial x_m}\right) (.,y)dy\ast G_{\nu},
\end{array}
\end{equation}
where we use Einstein notation for the first order spatial derivative of the kernel,
 Define
\begin{equation}
\delta v^{t_0,(0)}_i(t,x):=0,~\delta v^{t_0,(k)}_i:=v^{t_0,(k)}_i-v^{t_0,(k-1)}_i~ \mbox{for $k\geq 1$.} 
\end{equation}
You may still prove  that for local time $t\in [t_0,t_1]$ and $m\geq 2$
\begin{equation}
\delta v^{t_0,(k)}(t,.)\in {\cal C}^{m(n+1)}_{pol,m},~v^{t_0,(k)}(t,.)\in {\cal C}^{m(n+1)-\mu}_{pol,m}
\end{equation}
for some $\mu\in (0,1)$, i.e., the order of polynomial decay, characterized by the function space ${\cal C}^{m(n+1)}_{pol,m}$, is inherited by $\sum_{k}\delta v^{t_0,(k)}(t,.)$ (where local time contraction implies independence of the upper bound of the iteration index $k$). This observation can be used to designing global schemes using damping effects of the convoluted initial data terms, external control or auto-control.
\end{rem}

\begin{rem} 
In this article we add the additional assumption of strong polynomial decay in order to ensure strong spatial polynomial decay of local solutions. In a former version of this paper we omitted this additional assumption, and it may be optional indeed. In any case the additional assumption simplifies the proof. Here, note that strong spatial polynomial decay is inherited in local time fixed point iteration scheme
(such that terms $K^*_{n,i}(s,.-z)w^{E}_{m,j}(s,z)w^E_{j,m}(s,z)$) and the Burgers term have strong polynomial decay for the data chosen). Furthermore, note that the latter kernel term compensates growth of the rational function $\frac{\Pi_{i=1}^n(1+x_i^2)}{(1+x_j^2)(1+x_m^2)}$ in the crucial Leray projection term
 \begin{equation}
\begin{array}{ll}
+\int_{K^s_{\rho}}\sum_{j,m=1}^n\frac{(\rho-t)\sqrt{\rho^2-t^2}^3}{\rho^2(1+x_j^2)(1+x_m^2)}\times\\
\\
\times\left( K_{n,i}^*(.-z)\right) \left( \frac{\partial w^E_m}{\partial y_j}\frac{\partial w^E_j}{\partial y_m}\right) (s,z)\frac{\Pi_{i=1}^n(1+x_i^2)}{(\rho-t)^n}dz.
\end{array}
\end{equation}
\end{rem} 
The iteration scheme in (\ref{Navleraysol}) and the incompressibility condition imply that on a local 
time interval $[t_0,t_0+\Delta],~ t_0\geq 0$ for $\Delta >0$ we hav the representation 
\begin{equation}\label{Navlerayscheme2iaber}
\begin{array}{ll}
 v^{\nu}_i=v^{\nu}_i(t_0,.)\ast_{sp}G_{\nu}
-\sum_{j=1}^D \left( v^{\nu}_jv^{\nu}_i\right) \ast G_{\nu,j}\\
\\+\left( \int_{{\mathbb R}^D}\left( K_D(.-y)\right) \sum_{j,m=1}^D\left( \frac{\partial v^{\nu}_m}{\partial x_j}\frac{\partial v^{\nu}_j}{\partial x_m}\right) (.,y)dy\right) \ast G_{\nu,i},
\end{array}
\end{equation}
We have observed elsewhere that this representation of the increment  $|v^{\nu}_i-v^{\nu}_i(t_0,.)\ast_{sp}G_{\nu}|$ in terms of the first order spatial dervatives of the Gaussian, i.e., in terms of  $G_{\nu,i}$, imply especially  that for regular data in ${\cal C}^{m(n+1)}_{pol,m}$ we have for some $\mu\in (0,1)$ and some Lipschitz constant $L$ of the Burgers and the Leray projection term function  in (\ref{Navlerayscheme2iaber}) imply that
\begin{equation}
|v^{\nu}_i(t_0+\Delta,0)-v^{\nu}_i(t_0,.)\ast_{sp}G_{\nu}|\leq L\Delta^{1+\mu},
\end{equation} 
such that $|v^{\nu}_i(t_0+\Delta,0)|\neq 0$ can be enssured on smalltime intervals $\Delta >0$ for data with $|v^{\nu}_i(t_0,0)|\neq 0$. This is especially easy to see for symmetric data (but it holds in general).  
Therefore, in the following theorem the mentioned symmetric data refer to functions $h_i$ which satisfy
\begin{equation}\label{sym}
h_i(y)=h_i(y^{-j})~\mbox{for all $1\leq i,j\leq n$, cf. notation at (\ref{Lip}) below. }
\end{equation} 
\begin{thm}\label{lemeul}
Let $m\geq 2$ be an integer number and let $\epsilon >0$ be any small positive real number. We assume $n\geq 3$. For some large data $h_i\in H^{m}\cap {\cal C}^{m(n+1)}_{pol,m}$ for $m>\frac{n}{2}+1$ (and especially for some large symmetric data in the sense of (\ref{sym})) ,  for all $1\leq i\leq n$ there is a $\rho >0$ such that there exists a regular solution $w^E_i\in C^1\left(\left(0,\infty\right),H^r\right),~1\leq i\leq n$ of the equation (\ref{eulertrans}) for $r=m-\epsilon$ with initial data $h^{\rho}_i,~1\leq i\leq n$. Furthermore, there are data $h^{\rho}_i,~1\leq i\leq n$ corresponding to data $h_i\in H^{m}\cap {\cal C}^{m(n+1)}_{pol,m}$ for $m>\frac{n}{2}+1$ in original coordinates, such that  $h^{\rho}_i(0)\neq 0$ for some $1\leq i\leq n$, and such that there is a $\rho>0$ such that for small $\epsilon >0$ there is a solution $w^E_i\in C^1\left(\left(0,\infty\right),H^r\right),~1\leq i\leq n,~r=m-\epsilon$ of the associated Cauchy problem in (\ref{eulertrans}) which satisfies $\lim_{s\uparrow \infty}w^E_i(s,0)\neq 0$ for some $1\leq i\leq n$ at time $t(s)=\rho>0$ (corresponding to $s=\infty$). Hence,  the corresponding local solution function $v^E_i,~1\leq i\leq n$ is a local classical solution of the incompressible Euler equation and has a singularity at the point $(\rho,0)$, in the sense that
\begin{equation}
t\uparrow  \rho~\Rightarrow 
~{\Big |}v^E_i(t,0){\Big |}\uparrow \infty,
\end{equation}
where the singularity is at most of order $-1$, i.e.,
\begin{equation}
\sup_{0\leq t\leq \rho}{\Big |}v^E_i(t,0)(\rho-t){\Big |}\leq C<\infty.
\end{equation}
\end{thm}
\begin{rem}
We note that the limit $\lim_{s\uparrow \infty}w^E_i(s,0)$ exists for all $1\leq i\leq n$ and corresponds to the limit $\lim_{t\uparrow \rho}w^{*,E}_i(t,0)$, where
\begin{equation}
w^{*,E}_i(t,.)=w^E_i(s,.)
\end{equation}
for all $1\leq i\leq n$ with $s=\frac{t}{\sqrt{\rho^2-t^2}}\in [0,\infty)$ and $t\in [0,\rho)$. Due to existing limits the function $w^{*,E}_i(t,.),~1\leq i\leq n$ can be extended to the domain with time interval $[0,\rho]$.
\end{rem}

\begin{proof}
 A direct local iteration scheme seems to be inappropriate in order to prove contraction and existence, since the order of spatial derivatives in the representations of iterative approximations  increases at each iteration step for such schemes. Therefore we add a Laplacian operator $-\nu\Delta$  (times viscosity $\nu>0$) to the left side of the equation in (\ref{eulertrans}) and consider theviscosity limit $\nu\downarrow 0$ in a further step. We proceed in three steps. First in a) we prove existence for the family of equations (a family with parameter $\nu >0$)
\begin{equation}\label{eulertransnu}
\begin{array}{ll}
w^{\nu,E}_{i,s}-\nu \Delta w^{\nu,E}_i+\sum_{j=1}^n\frac{\sqrt{\rho^2-t^2}^3}{\rho^2(1+x_j^2)}w^{\nu,E}_j\frac{\partial w^{\nu,E}_i}{\partial y_j}\\
\\
-\sum_{j=1}^n\frac{(\rho-t)\sqrt{\rho^2-t^2}^3}{\rho^2}\arctan(x_j)w^{\nu,E}_{,j}(s,y)
=-\frac{\sqrt{\rho^2-t^2}^3}{\rho^2(\rho-t)}w^{\nu,E}_i\\
\\
+\sum_{j,m=1}^n\int\frac{(\rho-t)\sqrt{\rho^2-t^2}^3}{\rho^2(1+x_j^2)(1+x_m^2)}\left( K^*_{n,i}(.-z)\right) \left( \frac{\partial w^{\nu,E}_m}{\partial y_j}\frac{\partial w^{\nu,E}_j}{\partial y_m}\right) (s,z)\frac{\Pi_{i=1}^n(1+x_i^2)}{(\rho-t)^n}dz,
\end{array}
\end{equation}
for strong initial initial data $w^{\nu,E}_i(0,.)=h^{\rho}_i(.),~1\leq i\leq n$.
 More precisely we prove existence for closely related problems defined on the whole space (cf. item a) below). In a second step b) we prove that there exists a $\rho >0$ and data $h^{\rho}_i,~1\leq i\leq n$ with $w^{\nu,E}_i(0,0)=\rho h_i(0)= h^{\rho}_i(0)\neq 0$ for some $1\leq i\leq n$ such that
\begin{equation}\label{w*}
~w^{*,\nu,E}_i(\rho,0)\neq 0,
\end{equation}
where we recall that $w^{*,\nu,E}_i(t,.):=w^{\nu,E}_i(s,.)$. This may be achieved for a considerable class of data due to upper bound estimates for convolutions of first order derivatives of the Gaussian convoluted with Lipschitz continuous nonlinear terms  in the local fixed proint iteration scheme. However, for the proof of the weaker statement of Theorem \ref{lemeul} above it is sufficient to consider a specific data. For example for $\rho >0$ small choose
\begin{equation}\label{expdata}
h_i=\frac{1}{\rho}\exp\left(-\frac{\sum_{i=1}^3x_i^2}{a} \right),~a=\frac{1}{\rho^3} 
\end{equation} 
Then we have $w^{\nu,E}_i(0,0)=\rho h_i(0)= h^{\rho}_i(0)=1$, but for the first order derivatives 
\begin{equation}\label{upprho}
\max_{1\leq i,j\leq n}\sup_{x\in {\mathbb R}^n}|h_{i,j}(x){|}\lesssim \sqrt{\rho},
\end{equation}
which becomes small for small $\rho$. 
\begin{rem}
Note that the first order spatial derivative of $h_i$ in (\ref{expdata}) is 
\begin{equation}
|h_{i,j}(x){|}=\frac{1}{\rho}{\Bigg |}\frac{-2x_j}{a}\exp\left(-\frac{\sum_{i=1}^3x_i^2}{a} \right){\Bigg |}
\end{equation}
such that the maximal supremum $\max_{1\leq i,j\leq n}\sup_{x\in {\mathbb R}^n}|h_{i,j}(x){|}$ is obtained for some zero of this derivatives, i.e.,
\begin{equation}
\frac{1}{\rho}\left( \frac{-2}{a}+\frac{4x_j^2}{a^2}\right) \exp\left(-\frac{\sum_{i=1}^3x_i^2}{a} \right)=0~\mbox{or}~x_j=\pm \sqrt{\frac{a}{2}}.
\end{equation} 
Hence $|h_{i,j}(x){|}\leq \frac{1}{\rho}{\Bigg |}\frac{\sqrt{2}}{\sqrt{a}} {\Bigg |}$ for all $1\leq i,j\leq n$ and $x\in {\mathbb R}^n$ such that the upper bound in (\ref{upprho}) holds for  the choice of $a$ in (\ref{expdata}). 
\end{rem}
For such data the Burgers term and the Leray projection term of the transformed equation for $w^{\nu,E}_i$ become small and are easily estimated in the fixed point iteration scheme such that the solution increment over the time interval $[0,\rho]$ is small compared to $1$ and  (\ref{w*}) holds. 
In a third  step c) we prove that the properties of the solution $w^{\nu,E}_i,~1\leq i\leq n$ described in a) and  b), i.e., existence of a solution of the Cauchy problem in (\ref{eulertransnu}) and the property in (\ref{w*}) are preserved in the viscosity limit $\nu \downarrow 0$. Finally, we  conclude that the corresponding solution of the incompressible Euler equation has a singularity at $(\rho,0)$.
\begin{itemize} 

\item[a)] We define a fixed point iteration scheme based on the equation in (\ref{eulertransnu}). Recall that in (\ref{eulertransnu}) we understand $x_i=x_i(s,y_i)$ as functions of the transformed variables. Especially, in the Leray projection term $x_i$ is part of the convolution, i.e., the Leray projection term reads
\begin{equation}\label{lerayfirst}
\begin{array}{ll}
\sum_{j,m=1}^n\int\frac{(\rho-t)\sqrt{\rho^2-t^2}^3}{\rho^2(1+x^2_j(s,z_j))(1+x_m^2(s,z_m))}\left( K^*_{n,i}(.-z)\right) \left( \frac{\partial w^{\nu,E}_m}{\partial y_j}\frac{\partial w^{\nu,E}_j}{\partial y_m}\right) (s,z)\times\\
\\
\times \frac{\Pi_{i=1}^n(1+x_i^2(s,z_i))}{(\rho-t)^n}dz,
\end{array}
\end{equation}
where $x^2_i(s,z_i)$ denote the squared value of $x_i(s,z_i)$. Note that the integral in
(\ref{lerayfirst}) is with respect to the time section $K^s_{\rho}$ of the cone $K_{\rho}$ at each time $s$. We suppress this reference of the integral in the following for simplicity of notation (if the reference is clear from the context). We remark that at each time $s$ the spatial domain of the integral is 
\begin{equation}
-(\rho-t)\frac{\pi}{2}\leq z_i\leq (\rho-t)\frac{\pi}{2},~1\leq i\leq n,
\end{equation}
such that $dz=(\rho -t)^ndz^*$ for $z^*_i=\frac{z_i}{\rho-t}$. Hence, the factor $\left( \frac{1}{\rho-t}\right) ^n$ cancels if we consider a spatial transformation to a cube at each time.
Furthermore, the functions $w^{\nu,E}_i,~1\leq i\leq n$ and their derivatives have a strong spatial decay such that the spatial factor 
\begin{equation}
\Pi_{i=1}^n(1+x_i^2(s,z_i))
\end{equation}
in (\ref{lerayfirst}) is compensated by the convolution of the Laplacian kernel with these functions $w^{\nu,E}_i,~1\leq i\leq n$, and by the factor $\frac{1}{\rho^2(1+x^2_j(s,z_j))(1+x_m^2(s,z_m))}$ in (\ref{lerayfirst}). 
Note that the functions $w^{E}_i$ are supported on the cone $K_{\rho}$. We may extend $w^{E}_i$ trivially assuming that
\begin{equation}
w^{E}_i(s,y)=0\mbox{ for }(s,y)\in [0,\infty)\times {\mathbb R}^n\setminus K_{\rho},
\end{equation}
where we use the same symbol for these trivial extensions. The reason for this extension is that  approximations $w^{\nu,E}_i$ of the functions $w^{E}_i$ based on convolutions with the Gaussian are naturally defined on the whole space. They have exponential decay outside the cone and in the limit $\nu\downarrow 0$ the support is the cone $K_{\rho}$ of course. 
We prove that there exist solutions $w^{\nu ,E}_i,~1\leq i\leq n,~\nu>0$ of (\ref{eulertransnu}) (as an equation on the whole domain) with strong data $h_i,~1\leq i\leq n$ at time $s=0$ such that
\begin{equation}\label{wnu}
w^{\nu ,E}_i\in C^{1}
\left([0,\infty),H^{m-\epsilon}\left({\mathbb R}^n\right)\right)\cap C^1,~m-\epsilon>\frac{n}{2}+1.
\end{equation}
Here, we understand $H^{m-\epsilon}=H^{m-\epsilon}\left({\mathbb R}^n\right)$.
We work with classical representations of solutions of the equation in (\ref{eulertransnu}) in terms of convolutions with the Gaussian.
Let $G_{\nu}$ be the fundamental solution of the equation
\begin{equation}
\frac{\partial G_{\nu}}{\partial t}-\nu \Delta G_{\nu}=0.
\end{equation}
We have
\begin{equation}
G_{\nu}(t,x;s,y)=\frac{1}{\sqrt{4\pi \nu(t-s)}^n}\exp\left(-\frac{(x-y)^2}{4\nu(t-s)}\right).
\end{equation}
We consider the Gaussian $G_{\nu}$ on the whole domain such that $G_{\nu}$ is defined for all $(t,x)\in (0,\infty)\times {\mathbb R}^n$ and $$(s,y)\in \left\lbrace (s',y')|(s',y')\in [0,\infty)\times {\mathbb R}^n~\&~s<t\right\rbrace .$$
We note again: in order to work with convolutions with Gaussians on the whole space while the convoluted functions are defined only on a cone we (trivially) extend the latter functions to the whole space, and work with approximations of solutions to (\ref{eulertransnu}). We may do this as we are interested in the viscosity limit $\nu\downarrow 0$, and this viscosity limit is supported on the cone $K_{\rho}$. 
In the following we denote spatial convolutions with $G_{\nu}$ by $\ast_{sp}$ and convolutions with respect to space and time  by $\ast$.
 An approximative classical solution representation of the equation in (\ref{eulertransnu}) (while abbreviating $t=t(s)$) for $0<s<\infty$ corresponding to $0<t<\rho$ is
\begin{equation}\label{eulertransnu0}
\begin{array}{ll}
w^{\nu,E,\epsilon}_{i}(s,.)=w^{\nu,E,\epsilon}_i(0,.)\ast_{sp} G_{\nu}-\sum_{j=1}^n\frac{\sqrt{\rho^2-t^2}^3}{\rho^2(1+x_j^2)}\left( w^{\nu,E,\epsilon}_j\frac{\partial w^{\nu,E,\epsilon}_i}{\partial y_j}\right)\ast G_{\nu}\\
\\
+\sum_{j=1}^n\frac{(\rho -t)\sqrt{\rho^2-t^2}^3}{\rho^2}\arctan(x_j)w^{\nu,E,\epsilon}_{,j}\ast G_{\nu}=-\frac{\sqrt{\rho^2-t^2}^3}{\rho^2(\rho-t)}w^{\nu,E,\epsilon}_i\ast G_{\nu}+\\
\\
\sum_{j,m=1}^n\int\frac{(\rho-t)\sqrt{\rho^2-t^2}^3}{\rho^2(1+x_j^2)(1+x_m^2)}\left( K^*_{n,i}(.-z)\right) \left( \frac{\partial w^{\nu,E,\epsilon}_m}{\partial y_j}\frac{\partial w^{\nu,E,\epsilon}_j}{\partial y_m}\right) (s,z)\times \\
\\
\times \frac{\Pi_{i=1}^n(1+x_i^2)}{(\rho-t)^n}dz\ast G_{\nu},~\mbox{with}~(s,.)\in [0,\infty)\times {\mathbb R}^n,
\end{array}
\end{equation}
and where $w^{\nu,E,\epsilon}_i(0,.)=w^{\nu,E}_i(0,.),~1\leq i\leq n$. The parameter $\epsilon >0$ is a parameter which reminds us that we are working with approximative solutions defined on the whole domain with exponential decay outside the cone. Here, the support outside the cone vanishes as $\epsilon,\nu\downarrow 0$ . Furthermore, note that the integral of the Leray projection term is mainly over the cone section $K^s_{\rho}$ at given time $s$ such that an additional factor $(\rho-t)^n=(\rho-t(s))^n$ appears upon spatial transformation of the cone section to a cylinder. For small $\rho >0$ the solution function $w^{\nu,E,\epsilon}_i,~ 1\leq i\leq n$, i.e., the fixed point in (\ref{eulertransnu0}), can be constructed by local contraction in strong function spaces. The function $w^{\nu,E,\epsilon}_i,~1\leq i\leq n$ is an approximative solution of (\ref{eulertransnu}) as we have spatial exponential decay outside the cone $K_{\rho}$ which becomes stronger as the parameter $\nu$ goes to zero, such that the solution function values become identical to zero outside the cone as $\nu,\epsilon\downarrow 0$. 
For $m\geq \frac{1}{2}n+1$ and  for each $\nu >0$, $1\leq |\alpha|\leq m,~1\leq |\beta|+1= |\alpha|$ with $\beta_j+1=\alpha_j\neq 0$ and $\beta_l=\alpha_l$ for $l\neq j$ spatial derivatives of a regular solution $w^{\nu,E,\epsilon}_i,~1\leq i\leq n$ have the representation
\begin{equation}\label{eulertranslim}
\begin{array}{ll}
D^{\alpha}_y w^{\nu,E,\epsilon}_i=D^{\alpha}_yw^{\nu,E,\epsilon}_i(0,.)\ast_{sp} G_{\nu}\\
\\
+\sum_{j=1}^nD^{\beta}_y\left( \frac{\sqrt{\rho^2-t^2}^3}{\rho^2(1+x_j^2)} \left( w^{\nu,E,\epsilon}_j\frac{\partial w^{\nu,E,\epsilon}_i}{\partial y_j}\right)\right)  \ast G_{\nu,j}\\
\\
+\sum_{j=1}^n\frac{(\rho-t)\sqrt{\rho^2-t^2}^3}{\rho^2}D^{\beta}_y\left( \arctan(x_i) w^{\nu,E,\epsilon}_{,j}(s,y)\right) \ast G_{\nu,j}\\
\\
-\frac{\sqrt{\rho^2-t^2}^3}{\rho^2(\rho-t)}D^{\alpha}_x w^{\nu,E,\epsilon}_i\ast G_{\nu}
+\sum_{j,m=1}^n\int D^{\beta}_y{\Big (}\frac{(\rho-t)\sqrt{\rho^2-t^2}^3}{\rho^2(1+x_j^2)(1+x_m^2)}\times\\
\\
\times\left( K^{*}_{n,i}(y)\right)\left( \frac{\partial  w^{\nu,E,\epsilon}_m}{\partial y_j}\frac{\partial w^{\nu,E,\epsilon}_j}{\partial y_m}\right) (.,z)\frac{\Pi_{i=1}^n(1+x_i^2)}{(\rho-t)^n}dz{\big )}\ast G_{\nu,j}.
\end{array}
\end{equation}
Such classical representations (and the existence of a solution) can be justified by a contraction principle for an iteration scheme $\left( w^{\nu,E,\epsilon,k}_i\right)_{k\geq 0},~1\leq i\leq n$
of successive approximations of $w^{\nu,E,\epsilon}_i,~1\leq i\leq n$, which are defined recursively. At step $k=0$ we define
\begin{equation}
w^{\nu,E,\epsilon,0}_i=w^{\nu,E,0}_i=h^{\rho}_i\ast_{sp}G_{\nu},~1\leq i\leq n,
\end{equation}
and  for $k\geq 1$, the function $w^{\nu,E,\epsilon,k}_i,~1\leq i\leq n$ is defined recursively as the approximative (approximative solution in the sense of classical representations as outlined above) solution of the family of Cauchy problems
\begin{equation}\label{eulertransk}
\begin{array}{ll}
w^{\nu,E,\epsilon,k}_{i,s}+\sum_{j=1}^n\frac{\sqrt{\rho^2-t^2}^3}{\rho^2(1+x_j^2)}w^{\nu,E,\epsilon,k-1}_j\frac{\partial w^{\nu,E,\epsilon,k-1}_i}{\partial y_j}-\nu \Delta w^{\nu,E,\epsilon,k}_i\\
\\
+\sum_{j=1}^n\frac{(\rho-t)\sqrt{\rho^2-t^2}^3}{\rho^2}\arctan(x_j)w^{\nu,E,\epsilon,k-1}_{,j}(s,y)=-\frac{\sqrt{\rho^2-t^2}^3}{\rho^2(\rho-t)}w^{\nu,E,\epsilon,k}_i\\
\\
+\sum_{j,m=1}^n\int\frac{(\rho-t)\sqrt{\rho^2-t^2}^3}{\rho^2(1+x_j^2)(1+x_m^2)}\left( K_{n,i}(x-y)\right) \times\\
\\
\times\left( \frac{\partial w^{\nu,E,\epsilon,k-1}_m}{\partial y_j}\frac{\partial w^{\nu,E,\epsilon,k-1}_j}{\partial y_m}\right) (s,z)\frac{\Pi_{i=1}^n(1+x_i^2)}{(\rho-t)^n}dz.
\end{array}
\end{equation} 
For $k\geq 1$ and $1\leq i\leq n$ we define the functional increments
\begin{equation}
\delta w^{\nu,E,\epsilon,k}_i:=w^{\nu,E,\epsilon,k}_i-w^{\nu,E,\epsilon,k-1}_i,
\end{equation}
where we have 
\begin{equation}
\delta w^{\nu,E,\epsilon,k}_i(0,.)\equiv 0,
\end{equation}
and
\begin{equation}\label{eulertranslk*}
\begin{array}{ll}
\delta w^{\nu,E,\epsilon,k}_i=-\sum_{j=1}^n\frac{\sqrt{\rho^2-t^2}^3}{\rho^2(1+x_j^2)} \left( \delta w^{\nu,E,\epsilon,k-1}_j\frac{\partial w^{\nu,E,\epsilon,k-1}_i}{\partial y_j}\right) \ast G_{\nu}\\
\\
-\sum_{j=1}^n \frac{\sqrt{\rho^2-t^2}^3}{\rho^2(1+x_j)^2} \left( w^{\nu,E,\epsilon,k-1}_j\frac{\partial \delta w^{\nu,E,\epsilon,k-1}_i}{\partial y_j}\right) \ast G_{\nu}-\frac{\sqrt{\rho^2-t^2}^3}{\rho^2(\rho-t)}\delta w^{\nu,E,\epsilon,k}_i\ast G_{\nu}\\
\\
+\sum_{j=1}^n\frac{\sqrt{\rho^2-t^2}^3}{\rho^2}y_j\delta w^{\nu,E,\epsilon,k-1}_{,j}(s,y)\ast G_{\nu}
+2\int\sum_{j,m=1}^n\frac{(\rho-t)\sqrt{\rho^2-t^2}^3}{\rho^2(1+x_j^2)(1+x_m^2)}\times\\
\\
\times\left( K_{n,i}(.-z)\right) \left( \frac{\partial \delta w^{\nu,E,\epsilon,k-1}_m}{\partial y_j}\frac{\partial w^{\nu,E,\epsilon,k-1}_j}{\partial y_m}\right) (.,z)\frac{\Pi_{i=1}^n(1+x_i^2)}{(\rho-t)^n}dz\ast G_{\nu},
\end{array}
\end{equation}
and where the spatial integral of the Leray projection term has its main mass in the cone section $K^s_{\rho}$ at each time $s$ (up to small $\epsilon$ related to small $\nu$). For the derivatives of order $\alpha\neq 0$ we get analogous representations ( which are easily derived from (\ref{eulertranslim})). There are two linear terms 
in the representation (\ref{eulertranslk*}) and similar representations for spatial derivatives of order $\alpha$. One is the damping term of the form
\begin{equation}
-\frac{\sqrt{\rho^2-t^2}^3}{\rho^2(\rho-t)}D^{\alpha}_x\delta w^{\nu,E,\epsilon,k}_i\ast G_{\nu}
\end{equation}
which lowers the value function and its derivatives pointwise. 
For $\alpha=0$ we have
\begin{equation}
\begin{array}{ll}
\lim_{\nu\downarrow 0}\int_{K_{\rho}}\frac{\sqrt{\rho^2-t^2}^3}{\rho^2(\rho-t)}
|w^{\nu,E}_i(s,y)|ds dy\\
\\
=\lim_{\nu\downarrow 0}\int_{Z}\frac{\sqrt{\rho^2-t^2}^3}{\rho^2(\rho-t)}
|w^{*,\nu,E}_i(t,.)|\frac{\rho^2}{\sqrt{\rho^2-t^2}^3}(\rho-t)^ndt dz\\
\\ =\lim_{\nu\downarrow 0}\int_{Z}\frac{(\rho-t)^n}{(\rho-t)}
|w^{*,\nu,E}_i(t,.)|dt dz\leq C\mbox{vol}(Z)\frac{\left( \rho\right)^n}{n} \downarrow 0\\
\end{array}
\end{equation}
as $\rho \downarrow 0$, where $C>0$ is an upper bound of $\sup_{0\leq t\leq \rho}|w^{*\nu,E}_i(t,.)|$ and $\mbox{vol}(Z)$ is the volume of the cylinder 
$$Z=\left\lbrace (t,z)|(s,y)\in K_{\rho} \right\rbrace .$$
Note that $t\in [0,\rho]$ and $z\in \left(-\frac{\pi}{2},\frac{\pi}{2} \right)$ for elements of the cylinder. Here
\begin{equation}\label{note}
\frac{ds}{dt}=\frac{\rho^2}{\sqrt{\rho^2-t^2}^3},~z_i=\frac{y_i}{(\rho-t)}=\arctan(x_i),~1\leq i\leq n,
\end{equation}
such that
\begin{equation}
dy=\Pi_{i=1}^ndy_i=\Pi_{i=1}^n\frac{dy_i}{dz_i}dz_i=(\rho-t)^ndz.
\end{equation}
Hence, for small $\rho$ the damping is not strong enough such it can force that $w^{*,\nu,E}_i(\rho,0)=0$. The change of measure and the cylinder are considered above. 

The other linear term in (\ref{eulertranslk*}) is 
\begin{equation}
\sum_{j=1}^n\frac{\sqrt{\rho^2-t^2}^3}{\rho^2}y_j\delta w^{\nu,E,\epsilon,k-1}_{,j}(s,y)\ast G_{\nu}.
\end{equation}
Here we note that we understand $y_j=(\rho-t)\arctan(x_j)$ is bounded and can be kept bounded as we interpret $\delta w^{\nu,E,\epsilon,k-1}_{,j}$ to be defined on the whole domain. 
Exponential spatial decay is preserved at each iteration step $k$ of the scheme by the functional increments which have almost all their mass on the cone for data $h^{\rho}_i$. 
Note that nonlinear terms preserve exponential spatial decay a fortiori as products of functions with exponential spatial decay have stronger exponential spatial decay than their factors and this additional decay is stronger as the lowering effect of exponential decay caused by convolutions with the Gaussian or by convolutions with first order derivatives of the Gaussian.  In the Leray projection term the factor $\Pi_{i=1}^n(1+x_i^2)$ is compensated by the exponential spatial decay of the convolution $\left( K^*_{n,i}(.-z)\right) \left( \frac{\partial \delta w^{\nu,E,\epsilon,k-1}_m}{\partial y_j}\frac{\partial w^{\nu,E,\epsilon,k-1}_j}{\partial y_m}\right) (.,z)$ (and also by the factor  $\frac{1}{(1+x_j^2)(1+x_m^2)}$). 
Here recall that the kernel $K^*_{n,i}$ defined in (\ref{k*}) has a factor $\frac{x_i(y_i)}{|x(y)||^n}$.  We have
we have for all $k\geq 1$ and $s\geq 0$
\begin{equation}
\delta w^{\nu,E,\epsilon,k}_i(s,.)\in {\cal C}^{m(n+1)}_{pol,m}.
\end{equation}
For the Leray projection term we remark that a transformation of the cone $K^s_{\rho}$ to the cylinder $Z$ at each time $s$ ensures that the factor $\frac{1}{(\rho-t)^n}$ cancels such that the integrated Leray projection term can have an upper bound. 
The convolutions with the Gaussian, the first order spatial derivatives of the Gaussian, and the Laplacian kernel can be estimated by Young inequalities using the classical representations for $\delta w^{\nu,E,\epsilon,k}_i(s,.)$. The integrals over time $s$ are splitted into local time integrals for $s\in [0,1]$ and global time integrals for $s\geq 1$. The local time integrals ($s\in [0,1]$ are splitted into local spatial integrals with factor $1_{B_1}$ and their complements. For the local time and local spatial intergals we may use local standard estimates for the Gaussian and first order derivatives of the Gaussian, which are locally $L^1$. 
First for $\nu >0$ we have for $t\neq s$ and $x\neq y$ and $\mu \in (0,1)$ an upper bound for $G_{\nu}$ of the form
\begin{equation}\label{G1}
\begin{array}{ll}
{\Big |}\frac{1}{\sqrt{4\pi \nu(t-s)}^n}\left(\frac{(x-y)}{2\sqrt{\nu(t-s)}}\right)^{2\mu-n}\left(\frac{(x-y)}{2\sqrt{\nu(t-s)}}\right)^{n-2\mu} \exp\left(-\frac{(x-y)^2}{4\nu(t-s)}\right){\Big |}\\
\\
\leq \frac{C}{\sqrt{\nu(t-s)}^{\mu}|x-y|^{n-2\mu}},~ \mbox{where}~
C=\frac{1}{\sqrt{\pi}^n}\sup_{z>0}z^{n-2\mu}\exp(-z^2).
\end{array}
\end{equation}
Note that the upper bound in (\ref{G1}) is integrable for $\mu\in (0,1)$. For the first spatial derivatives $G_{,i}$ we have an additional factor which we may estimate by $
{\big |}\frac{(x_i-y_i)}{\mu(t-s)}{\big |}\leq \frac{1}{{|}x_i-y_i{|}}\frac{(x-y)^2}{\mu(t-s)}$
such that the upper bound of $G_{,i}$ of the form
\begin{equation}
\frac{C'}{\sqrt{\nu(t-s)}^{\mu}|x-y|^{n+1-2\mu}}
\end{equation}
holds for $t\neq s$ and $x\neq y$ and becomes locally integrable for $\mu \in (0.5,1)$.
The functions $D^{\alpha}_xw^{\nu,E,\epsilon,k}_i(s,.)$ are $L^2\cap C$ for all $0\leq |\alpha|\leq m$ inductively for all $k\geq 0$ where products of functions can be estimated by suprema of one factor (we do not even need the product rule for regular Sobolev norms). For the complementary time local ($s\in [0,1]$) and spatially global estimates we may use the exponential decay of the truncated Gaussian
\begin{equation}\label{G11}
\begin{array}{ll}
{\Big |}1_{{\mathbb R}^n\setminus B_1}\frac{1}{\sqrt{4\pi \nu(t-s)}^n}\exp\left(-\frac{(x-y)^2}{4\nu(t-s)}\right){\Big |}.
\end{array}
\end{equation}
Here in any case $1_{B_1}$ is the characteristic function which is $1$ on the ball of radius $1$ and zero elsewhere, and $1_{{\mathbb R}^n\setminus B_1}$ is the complementary characteristic function.
For the leray projection term we note that the truncated kernel $1_{B_1}K_{n,i}$  is in $L^1$ while the complement $1_{{\mathbb R}^n\setminus B_1}K_{n,i}$ is in $L^2$, and  the Gaussian estimates above , Plancherels' identity and Young inequalities ensure boundedness of the Leray projection term in local time.  We have given the details of this local time part of contraction results elsewhere, and need not to repeat every detail here. 
For large time estimates ($s\in ( 1,\infty)$) we may use a simple Gaussian upper bound based on the factor $\sqrt{s}^{-n}$, i.e.,
\begin{equation}\label{G11}
\begin{array}{ll}
{\Big |}\frac{1}{\sqrt{4\pi \nu s}^n}\exp\left(-\frac{(x-y)^2}{4\nu s}\right){\Big |}\leq \frac{1}{\sqrt{4\pi \nu s}^n},~ s\in (1,\infty),
\end{array}
\end{equation}
where the right side is global integrable for $n\geq 3$. Summing up we can extract a factor $\rho$ from all terms on the right side of (\ref{eulertranslk*}) such that
\begin{equation}\label{cstar}
\max_{1\leq i\leq n}{\big |}\delta w^{\nu,E,\epsilon,k}_i{\big |}_{H^m\cap C^m}\leq \rho c^*\max_{1\leq i\leq n}{\big |}\delta w^{\nu,E,\epsilon,k-1}_i{\big |}_{H^m\cap C^m},
\end{equation}
where for $\rho >0$ small we choose
\begin{equation}
c^*:=4\cdot2^m C_{h\rho}^2C_G(1+C_{K_n}),
\end{equation}
along with
\begin{equation}
{\big |}h^{\rho}_i{\big |}_{H^m\cap C^m}\leq C_{h\rho},
\end{equation}
and 
\begin{equation}
C_{K_n}={\big |}K_{,i}{\big |}_{L^1(B_1)}+{\big |}K_{,i}{\big |}_{L^2({\mathbb R}^n\setminus B_1)},
\end{equation}
and 
\begin{equation}
C=\max\left\lbrace C,C'\right\rbrace .
\end{equation}
The factor $2^m$ is due to the number of terms in the expansion of derivatives.
Note that we can choose
\begin{equation}
\max_{1\leq i\leq n}{\big |}h^{\rho}_i{\big |}_{H^m\cap C^m}\mbox{ close to $1$}
\end{equation}
for transformed data for convenience.

\item[b)] We consider a family of data functions $h^{\rho}_i,~1\leq i\leq n,~\rho >0$ such that for some $1\leq i_0\leq n$
 \begin{equation}
 \rho h_{i_0}(0)= h^{\rho}_{i_0}(0)=1,~h_i\in C^{m}\cap H^m\cap  {\cal C}^{m(n+1)}_{pol,m},~1\leq i\leq n,
 \end{equation}
where $h_i$ and $h^{\rho}_i,~1\leq i\leq n$ are the data in the respective coordinates and $\mbox{supp}(h^{\rho}_i)\subseteq  \left]-\rho\frac{\pi}{2},\rho\frac{\pi}{2}\right[^n$. For convenience we may choose data as in (\ref{expdata}).
\begin{rem}
There are two arguments for the fact that the functional increment $\delta w^{\nu,R,\epsilon}_i=w^{\nu,E,\epsilon}_{i}(s,.)-w^{\nu,E,\epsilon}_i(0,.)\ast_{sp} G_{\nu}$ becomes small as $\rho$ becomes small. In order to onserve this consider the terms of these functional increments in the representation in (\ref{eulertransnu0}), i.e., we have to consider the terms on the right side of 
\begin{equation}\label{eulertransnu0*}
\begin{array}{ll}
\delta w^{\nu,E,\epsilon}_i=-\sum_{j=1}^n\frac{\sqrt{\rho^2-t^2}^3}{\rho^2(1+x_j^2)}\left( w^{\nu,E,\epsilon}_j\frac{\partial w^{\nu,E,\epsilon}_i}{\partial y_j}\right)\ast G_{\nu}\\
\\
+\sum_{j=1}^n\frac{(\rho -t)\sqrt{\rho^2-t^2}^3}{\rho^2}\arctan(x_j)w^{\nu,E,\epsilon}_{,j}\ast G_{\nu}=-\frac{\sqrt{\rho^2-t^2}^3}{\rho^2(\rho-t)}w^{\nu,E,\epsilon}_i\ast G_{\nu}+\\
\\
\sum_{j,m=1}^n\int\frac{(\rho-t)\sqrt{\rho^2-t^2}^3}{\rho^2(1+x_j^2)(1+x_m^2)}\left( K^*_{n,i}(.-z)\right) \left( \frac{\partial w^{\nu,E,\epsilon}_m}{\partial y_j}\frac{\partial w^{\nu,E,\epsilon}_j}{\partial y_m}\right) (s,z)\times \\
\\
\times \frac{\Pi_{i=1}^n(1+x_i^2)}{(\rho-t)^n}dz\ast G_{\nu},~\mbox{with}~(s,.)\in [0,\infty)\times {\mathbb R}^n.
\end{array}
\end{equation}
The simplified argument stasts the iteration scheme with data as in (\ref{expdata}). Note that in any case we have $w^{\nu,E,\epsilon}_i$ is close to $(\rho-t)v^{E}_i$ for small $\nu>0$, and for the choice of data in  (\ref{expdata}) we first observe that $\frac{\partial w^{\nu,E,\epsilon,0}_i}{\partial y_j}$ is of order $\sqrt{\rho}$ and contraction analyis of the iteration scheme show that this behavior is inherited form iteration step to iteration step. Here recall that $v^{E}_{i,j}(t,x)=\frac{ w^{E}_{i,j}}{1+x_j^2}(s,y)$. It follows that the terms on the right side become small as $ \rho >0$ becomes small. Here note that the factor  $\frac{\sqrt{\rho^2-t^2}^3}{\rho^2}$ is cancelled if we transform from time $s$ to time $t$.  
We also mention a second stronger  argument which shows that the data have not to be chosen that specifically. We have shon elsewhere that Lipschitz continuous data $F$ (such as local time Burgers- and Leray projection solution terms in our case) which are convolted with first order spatial derivatives of the Gaussian have  
\begin{equation}\label{estest}
\begin{array}{ll}
{\Big |} F\ast G_{\nu,i}=\int F(t-\sigma,x-y)\left( \frac{-2y_i}{4\nu  \sigma}\right) G_{\nu}(\sigma,y)dyd\sigma {\Big |}\\
\\
={\Big |} \int \int_{y_i\geq 0} \left( F(t-\sigma,x-y)-F(t-\sigma,x-y^-)\right)\left( \frac{-2y_i}{4\nu  \sigma}\right) G_{\nu}(\sigma,y)dyd\sigma{\Big |}\\
\\
\leq {\Big |} L\int \int_{y_i\geq 0} \left( \frac{4y^2_i}{4\nu   \sigma}\right) G_{\nu}(\sigma,y)dyd\sigma {\Big |}=4 LM_2,
\end{array}
\end{equation} 
where $M_2$ is a finite second moment constant of the Gaussian which is is related for a small time step size $\Delta$ to upper bound as smallas
\begin{equation}\label{parafreeest}
 {\Big |} L\int_{0}^{\Delta}   \frac{1}{4} \frac{1}{\sqrt{4\pi}^D}\frac{1}{\sigma^{2.5}}\exp\left(-\frac{1}{4\sigma} \right)d\sigma {\Big |}
\end{equation}
which becomes small for small $\Delta >0$. This has been shown elsewhere.
\end{rem}
We show that there exists a $\rho >0$ such that for data $h^{\rho}_i,~1\leq i\leq n$ and $w^{\nu,E}_i(0,0)=\rho h^{\rho}_i(0)= 1$ for some $1\leq i_0\leq n$ we have
\begin{equation}
w^{*,\nu,E}_{i_0}(\rho,0)=\lim_{s\uparrow\infty} w^{\nu,E}_i(s,0)\neq 0.
\end{equation}
We have for all $1\leq i\leq n$ and $s\in [0,\infty)$
\begin{equation}
 w^{\nu,E,\epsilon}_i(s,.)=h^{\rho}_i(.)+ \sum_{l=1}^{\infty}\delta w^{\nu,E,\epsilon,l}_i(s,.),
\end{equation}
where the contraction result with respect to the $\sup_{s\geq 0}|.|_{H^m\cap C^m}$-norm (with $m\geq 2$) implies that this representation holds pointwise for the functions $w^{\nu,E,\epsilon}_i,~1\leq i\leq n$ and for spatial derivatives of these functions up to order $m$. 
Next, since for $k\geq 2$, it follows from \ref{cstar} that for some $m\geq 2$, all $s\in [0,\infty)$, for a given $\mu\in (0,1)$ and $\rho>0$ small enough we have 
\begin{equation}\label{rhoeps1}
\max_{1\leq i\leq n}{\big |}\delta w^{\nu,E,\epsilon,k}_i(s,.){\big |}_{H^m\cap C^m}\leq \rho^{\mu}\max_{1\leq i\leq n}{\big |}\delta w^{\nu,E,\epsilon,k-1}_i(s,.){\big |}_{H^m\cap C^m},
\end{equation}
where $\rho$ can be chosen such that $\rho^{\mu}\leq 0.5$, such that
\begin{equation}\label{rhoeps}
 \sum_{l=2}^{\infty}{\big |}\delta w^{\nu,E,\epsilon,l}_i(s,.){\big |}_{H^m\cap C^m}
\leq \frac{\rho^{\mu}}{1-\rho^{\mu}}
 {\big |} \delta w^{\nu,E,\epsilon,1}_i(s,.){\big |}_{H^m\cap C^m}.
\end{equation}
Next we determine an upper bound for $\sup_{s\geq 0}{\big |}\delta w^{\nu,E,\epsilon,1}_i(s,.){\big |}_{H^m\cap C^m}$. 
From (\ref{eulertransk}) and $\delta w^{\nu,E,\epsilon,1}_i:=w^{\nu,E,\epsilon,1}_i-w^{\nu,E,\epsilon,0}_i=w^{\nu,E,\epsilon,1}_i-h_i\ast_{sp}G_{\nu}$ we have
\begin{equation}\label{eulertranslk*1}
\begin{array}{ll}
\delta w^{\nu,E,\epsilon,1}_i=-\sum_{j=1}^n\frac{\sqrt{\rho^2-t^2}^3}{\rho^2(1+x_j^2)}  h^{\rho}_j\frac{\partial h^{\rho}_i}{\partial y_j}\ast G_{\nu}-\frac{\sqrt{\rho^2-t^2}^3}{\rho^2(\rho-t)} w^{\nu,E,\epsilon,1}_i\ast G_{\nu}\\
\\
+\sum_{j=1}^n\frac{\sqrt{\rho^2-t^2}^3}{\rho^2}y_jh^{\rho}_{,j}(s,y)\ast G_{\nu}
+\int\sum_{j,m=1}^n\frac{(\rho-t)\sqrt{\rho^2-t^2}^3}{\rho^2(1+x_j^2)(1+x_m^2)}\times\\
\\
\times\left( K_{n,i}(.-z)\right) \left( \frac{\partial h^{\rho}_m}{\partial y_j}\frac{\partial h^{\rho}_j}{\partial y_m}\right) (.,z)\frac{\Pi_{i=1}^n(1+x_i^2)}{(\rho-t)^n}dz\ast G_{\nu}.
\end{array}
\end{equation}
Hence, for $\mu\in (0,1)$ and $\rho$ small enough
\begin{equation}\label{euleko}
\begin{array}{ll}
\sup_{s\geq 0}{\big |}\delta w^{\nu,E,\epsilon,1}_i(s,.){\big |}_{H^m\cap C^m}\leq\rho c^*\leq \rho^{\mu},
\end{array}
\end{equation} 
with $c^*$ as in item a).
Hence, for $\rho$ small enough such that
\begin{equation}
\rho^{\mu}\left(1+\frac{\rho^{\mu}}{1-\rho^{\mu}}\right) \leq \frac{1}{2}h^{\rho}_{i_0}(0)=1\neq 0
\end{equation}
for some small $\epsilon >0$ (cf. item a)), we conclude that
\begin{equation}
w^{*,\nu,E,\epsilon}_i(\rho,.):=\lim_{s\uparrow \infty}w^{\nu,E,\epsilon}_i(s,.)\neq 0.
\end{equation}
\begin{rem}
Note that the conclusion in item b) can also be obtained easily from
\begin{equation}\label{wnu}
w^{*,\nu ,E}_i\in C^{1}
\left([0,\rho],H^m\cap C^m\right),~\mbox{for some $m\geq 2$}.
\end{equation}
However, using the contraction result explicitly we have an explicit upper bound for $\rho$. 
\end{rem}

\item[c)] First we mention that the contraction constant in item a) can be chosen independently of the viscosity $\nu$. The reason for this is that in the classical representations for $D^{\alpha}w^{\nu,E,\epsilon}_i,~1\leq i\leq n$ for $0\leq |\alpha|\leq m$ are of the form
\begin{equation}
H\ast G_{\nu}~\mbox{( for $|\alpha|=0$)},~F\ast G_{\nu,i} \mbox{ (for $|\alpha|>0$)},
\end{equation}
for regular functions $H$ and $F$. The convolutions $H\ast G_{\nu}$ have a natural upper bound such that for sequences of functions with strong spatial polynomial decay we may use a transformation and compactness arguments (related to Rellich's theorem) and obtain a regular limit with a small loss of regularity (cf. below for the use of strong spatial polynomial decay in this context). For the functionals $F\ast G_{\nu,i}$ we may use Lipschitz continuity of the functional $F(t,.)\in L^2\cap C$ (the convoluted  Burgers and Leray projection functionals and their spatial derivatives at each iteration step)
\begin{equation}\label{Lip}
{\big |}F(x-y)-F(x-y'){\big |}=:{\big |}F_x(y)-F_x(y'){\big |}\leq L|y-y'|
\end{equation}
for some constant $L>0$, where for $y^{-i}=\left( y^{-i}_1,\cdots, y^{-i}_n\right) $ with $y^{-i}_k=y_k$ for $k\neq i$ and $y^{-i}_i=-y_i$ we have
\begin{equation}
\begin{array}{ll}
{\big |}F\ast G_{\nu,i}{\big |}={\Big  |}\int F(x-y)\frac{2y_i}{4t}G_{\nu}(t,y)dtdy{\Big |}\\
\\
{\Big  |}\int_{y_i\geq 0}(F_x(y)-F_x(y^{-i})) \frac{2y_i}{4t}G_{\nu}(t,y)dtdy{\Big |}
\leq L{\Big  |}\frac{4y^2_i}{4t}G_{\nu}(t,y)dtdy{\Big |}
\end{array}
\end{equation}
which leads to natural $\nu$-indpendent estimates in item a). Next we prove c1) that the viscosity limit exists for the function $w^{\nu,E,\epsilon}_i(s,.),~1\leq i\leq n$ , i.e., for $m\geq 2$ and $h^{\rho}_i\in H^m\cap C^m,~1\leq i\leq n$ we prove 
\begin{equation}\label{sol1statend}
w^{*,E,\epsilon}_i:=\lim_{\nu\downarrow 0}w^{\nu ,E,\epsilon}_i\in C^{1}\left([0,\rho],H^{m-\epsilon}\left({\mathbb R}^n\right)\right), 
\end{equation}
and where for some $m>\frac{n}{2}+1$ we have for some finite constant
\begin{equation}\label{sol1astatend2}
\sup_{0\leq t\leq \rho}\max_{1\leq i\leq n}{\big |}w^{*,E,\epsilon}_i(t,.){\big |}_{H^r}\leq C.
\end{equation}
Note that this implies 
\begin{equation}
\begin{array}{ll}
\sup_{0\leq t\leq \rho}{\big |}w^{*,E,\epsilon}_i(t,.){\big |}_{C^1}=\sup_{0\leq s< \infty}{\big |}w^{E,\epsilon}_i(s,.){\big |}_{C^1}\leq \tilde{C}
\end{array}
\end{equation}
for some finite constant $\tilde{C}>0$ by the Sobolev lemma,
 and where (as usual) we denote  $|f|_{C^1}=\sum_{0\leq |\alpha|\leq 1}\sup_{x\in {\mathbb R}^n}|D^{\alpha}_xf(.)|$. Here recall that the dummy upper script $\epsilon>0$ reminds us that this solution is the viscosity limit of a problem which is defined on the whole domain. We shall conclude that a corresponding short time classical solution $v^E_i,~1\leq i\leq n$ of the incompressible Euler equation exists.
In item a) we have constructed regular upper bounds of $w^{\nu,E,\epsilon}_i~1\leq i\leq n$ which are independent of the viscosity $\nu$, i.e., for each $m\geq 2$, regular data $h_i,~1\leq i\leq n$ (as assumed) and corresponding regular data $h^{\rho}_i,~1\leq i\leq n$ we have a finite constant $C_m$ which is independent of $\nu$ such that
\begin{equation}
\max_{1\leq i\leq n}\sup_{s\geq 0}{\big |}w^{\nu,E,\epsilon}_{i}(s,.){\big |}_{H^m\cap C^m}\leq C_m,
\end{equation}
and where in addition $w^{\nu,E,\epsilon}_i(s,.)\in {\cal C}^{m(n+1)}_{pol,m},¸1\leq i\leq n,~\mbox{for}~s\geq 0$.

The existence of a viscosity limit is based on compactness arguments and the Lipschitz-continuity of the  terms in (\ref{eulertransnu0}), which are convoluted with the Gaussian.

Here, we note that convergence in strong norms implies that compactness arguments are available. Indeed, we can apply Rellich's embedding as follows - note that Rellich's theorem is formulated on bounded domains. Recall that Rellich's embedding is with respect to spaces $H^s_0\left( \Omega\right) $ which are closures in $H^s$ for $s>0$ of $C^{\infty}_c(\Omega)$, i.e., the space of smooth function with compact support. These are the functions of $H^s$ which are supported in the closure of $\overline{\Omega}$ of $\Omega$. For simplicity let is consider again the integer value $m\geq 2$.  For $0\leq |\alpha|\leq m$ we consider the functions
\begin{equation}
\left]-\frac{\pi}{2},\frac{\pi}{2}\right[^n \ni y\rightarrow \left( D^{\alpha}_xw^{\nu,E,\epsilon}_i\right)^*(s,z)=D^{\alpha}_xw^{\nu,E,\epsilon}_i\left( s,\tan \left( \frac{y}{\rho-s}\right) \right) ,
\end{equation}
 where $\tan(y)=(\tan(y_1),\cdots ,\tan(y_n))^T$. Since (\ref{wiepsilon}) holds , and
 \begin{equation}\label{alphaderc}
{\big |}\left( D^{\alpha}_xw^{\nu,E,\epsilon}_i\right) ^*(s,z){\big |}=c\left(1+|x|^{2m}\right) {\big |}D^{\alpha}_xw^{\nu,E,\epsilon}_{i}(s,x){\big |},
\end{equation} 
it follows that
 \begin{equation}\label{sol1stat*}
w^E_i:=\lim_{\nu,\epsilon\downarrow 0}w^{\nu ,E,\epsilon}_i\in C^{1}\left([0,\infty),H^r\left({\mathbb R}^n\right)\right) \mbox{ for }\frac{n}{2}<r<m.
\end{equation}
Note that for $m\geq  3$ we have $\nu \Delta w^{\nu,E,\epsilon}_i\downarrow 0$ as $\nu\downarrow 0$ while $w^{\nu,E}_i,~1\leq i\leq n$ is a classical solution of the dampened incompressible Navier Stokes equation with viscosity $\nu >0$. It follows that the viscosity limit $w^E_i,~1\leq i\leq n$ is a regular solution of a dampened incompressible Euler type equation corresponding to a local regular solution $v^{E}_i,~1\leq i\leq n$ of the original Euler equation on the time interval $[0,\rho)$. Together with local contraction this leads to
\begin{equation}
\sup_{0\leq t(s)\leq \rho}{\big |}w^E_i(s,.){\big |}_{H^r}:=\sup_{0\leq t(s)\leq \rho}\lim_{\nu\downarrow 0}{\big |}w^{\nu ,E}_i(s,.){\big |}_{H^r}\leq C
\end{equation}
for some constant $C>0$. 
Note that the existence of a classical solution of the incompressible Euler equation is implied. Indeed, iterating spatial derivatives starting from (\ref{firstder}) we get for all $t=t(s)<\rho$ and corresponding $s\geq 0$
\begin{equation}\label{alphaderc}
{\big |}D^{\alpha}_xv^E_{i}(t,.){\big |}=c\left(1+|x|^{2m}\right) {\big |}D^{\alpha}_xw^E_{i}(,.){\big |}
\end{equation} 
by induction, and local classical existence of a solution $v^E_i,~1\leq i\leq n$ on the time interval $[0,\rho)$ follows from compactness.
Finally c3) we show that the property
$w^{*,E}_{i_0}(\rho,0)\neq 0$ can be derived from 
\begin{equation}
\mbox{ for all }~\nu>0~w^{*,\nu,E}_{i_0}(\rho,0)\neq 0.
\end{equation}
This follows from the fact that in estimate in b) the finite constant $c^*$ can be chosen independently of the viscosity $\nu$.  In addition, we observe that
for $\rho$ small enough we have for some $\mu\in (0,1)$
\begin{equation}
\rho c^*\left(1+\frac{\rho^{\mu}}{1-\rho^{\mu}}\right) \leq \frac{1}{2}h_{i_0}(0)\neq 0
\end{equation}
for some small $\epsilon >0$ (cf. item a)).  We conclude that for he converging subsequence $w^{*,\nu_k,E,\epsilon}_{i_0}(\rho,.),~k\geq 1$ we have
\begin{equation}
w^{*,\nu_k,E,\epsilon}_{i_0}(\rho,.):=\lim_{s\uparrow \infty}w^{\nu_k,E,\epsilon}_{i_0}(s,.)\geq  \frac{h_{i_0}(0)}{2}.
\end{equation}
Hence, 
\begin{equation}
\lim_{k\uparrow \infty}w^{*,\nu_k,E,\epsilon}_{i_0}(\rho,.):=\lim_{s\uparrow \infty}w^{\nu_k,E,\epsilon}_{i_0}(s,.)\geq  \frac{h_{i_0}(0)}{2}\neq 0.
\end{equation}
It follows that $v^E_{i_0}(t,0)=\frac{w^E_i(\rho,0)}{\rho-t}$ has a singularity at $t=\rho$.
\end{itemize}
\end{proof}

The proof in Theorem \ref{lemeul} goes through for strong data $h_i\in H^{m}\cap C^m\cap {\cal C}^{m(n+1)}_{pol,m}$ for all $m\geq 2$, where the transformed data satisfy $ h^{\rho}_i(0)=1\neq 0$ for small $\rho$ and $w^{E}_{i_0}(\rho,0)\neq 0$. We conclude
\begin{cor}\label{coreul}
In the situation of the last statement of Theorem \ref{lemeul} there are data $h_i,~1\leq i\leq n$ such that $h_i\in C^{\infty}$ and such that the corresponding solution function $v^E_i$ of the incompressible Euler equation defined by $v^E(t,x)=\frac{1}{\rho-t}w^E_i(s,y),~1\leq i\leq n$ has a singularity at the point $(\rho,0)$.
\end{cor}
Next we consider extended models with time-dependent external force data $f_i\in L^2$ for $1\leq i\leq n$. Note that in any case such models include models with nonzero initial data implicitly, i.e., schemes with regular initial data $h_i\in H^s,~1\leq i\leq n$ for $s>\frac{n}{2}+1$ can essentially be always written in the form with zero initial data and an adjusted external force term. Just consider $v'_i=v_i-h_i$ instead of the velocity function $v_i$. This leads to additional linear first order derivative term in the equation for $v'_i,¸1\leq i\leq n$, but this is not essential.  Here 'essential' means that the additional terms in the equation for $v'_i,¸1\leq i\leq n$ do not alter the analysis of the equation essentially. Therefore, we may consider models with nonzero regular initial data in the following without loss of generality in the sense that the following considerations can also be applied to zero initial data and related external force terms in $L^2$. Recall from (\ref{firstder}) that the singularity factor $\frac{1}{\rho-t}$ in $v^E_i=\frac{w^E_i}{\rho -t}$ cancels for the first order spatial derivatives, i.e.,
\begin{equation}\label{firstder*}
v^E_{i,j}=\frac{w_{i,j}(s,y)}{1+x_j^2},
\end{equation}
hence we get ($x_j\equiv x_j(y_j)$)
\begin{equation}\label{vw}
\begin{array}{ll}
\Delta v^E_i=\sum_{j=1}^nv^E_{i,j,j}=\sum_{j=1}^n\frac{w_{i,j,j}(s,y)}{1+x_j^2}\frac{dy_j}{dx_j}-\frac{w_{i,j}(s,y)}{(1+x_j^2)^2}2x_j\\
\\
=\sum_{j=1}^n\frac{w_{i,j,j}(s,y)}{1+x_j^2}\frac{\rho-t}{1+x_j^2}-\frac{w_{i,j}(s,y)}{(1+x_j^2)^2}2x_j
\end{array}
\end{equation}
This means that $\Delta w^E_i\in L^2$ implies $\Delta v^E_i\in L^2$. Note that for data $h_{i}\in H^s$ with $s>\frac{n}{2}+2$ we even get solutions $w^E_i,~1\leq i\leq n$ with $w^E_i\in H^r$
for some order $\frac{n}{2}+2< r<s$, which implies that $w_{i,j,j}$ in (\ref{vw}) is of bounded modulus, which makes the latter conclusion even more obvious. In the following statement we refer to the latter observation as a 'situation of data'.
As a consequence of Lemma \ref{lemeul} and Corollary \ref{coreul} we get
\begin{thm}
In the situation of Theorem \ref{lemeul} or in a situation of data with higher regularity consider a solution function $w^E_i,~1\leq i\leq n$ corresponding to a solution $v^E_i,~1\leq i\leq n$ of the incompressible Euler equation as described in Corollary \ref{coreul}. Define
\begin{equation}\label{fw}
f^{w,E}_i=-\nu\Delta w^E_i.
\end{equation}
According to lemma \ref{lemeul} there exists a $\rho >0$ such that the function $w^E_i,~1\leq i\leq n$ satisfies (with $t=t(s)$)
\begin{equation}\label{eulertranslast}
\begin{array}{ll}
w^E_{i,s}-\nu\Delta w^E_i+\frac{\sqrt{\rho^2-t^2}^3}{\rho^2(1+x_j^2)}\sum_{j=1}^nw^E_j\frac{\partial w^E_i}{\partial y_j}\\
\\
-\sum_{j=1}^n\frac{(\rho-t)\sqrt{\rho^2-t^2}^3}{\rho^2}\arctan(x_j)w^{\nu,E}_{,j}(s,y)=-\frac{\sqrt{\rho^2-t^2}^3}{\rho^2(\rho-t)}w^E_i
+f^{w,E}_i(s,.)+\\
\\
\sum_{j,m=1}^n\int\frac{(\rho-t)\sqrt{\rho^2-t^2}^3}{\rho^2(1+x_j^2)(1+x_m^2)}\left( K^*_{n,i}(.-z)\right) \left( \frac{\partial w^{E}_m}{\partial y_j}\frac{\partial w^{E}_j}{\partial y_m}\right) (s,z)\frac{\Pi_{i=1}^n(1+x_i^2)}{(\rho-t)^n}dz,
\end{array}
\end{equation} 

with the initial data
\begin{equation}
w^E_i(0,.)=h_i,~1\leq i\leq n
\end{equation}
on the time interval $t\in[0,\rho)$ corresponding to the time interval $s\in [0,\infty)$, and where the solution is a solution in the classical sense. Furthermore the solution can be extended to $t=\rho$ at $y=0$, where $\lim_{s\uparrow \infty}w^E_i(s,0)=\lim_{t(s)\uparrow \rho}w^E_i(s,0)\neq 0$. The corresponding functions $v^E_i,~1\leq i\leq n$ satisfy
\begin{equation}\label{eulertransnavv}
\begin{array}{ll}
v^E_{i,t}-\nu\Delta v^E_i+\sum_{j=1}^nv^E_j\frac{\partial v^E_i}{\partial y_j}+\\
\\
=\sum_{j,m=1}^n\int_{{\mathbb R}^n}\left( \frac{\partial}{\partial x_i}K_n(x-y)\right) \sum_{j,m=1}^n\left( \frac{\partial v_m}{\partial x_j}\frac{\partial v_j}{\partial x_m}\right) (t,y)dy+f^{v,E}_i(t,.),
\end{array}
\end{equation} 
with the initial data
\begin{equation}
v^E_i(0,.)=h_i,~1\leq i\leq n,
\end{equation}
and where the relation of $f^{v,E}_i(t,.)$ and $f^{w,E}_i(s,.)$ is given via (\ref{fw}) and (\ref{vw}).
Here, we have that $f^{v,E}_i\in L^2$ and $v^E_i$ has a singularity at $(\rho,0)$ for an index $i$ with  $\lim_{s\uparrow \infty }w_i(s,0)\neq 0$.   
\end{thm}

\section{Conclusion}

The preceding considerations show that the argument  discussed in \cite{web} cannot lead to a unique solution in the spaces chosen there (even if the error outlined in \cite{web} can be corrected). For that purpose it is necessary that the external force terms are located in more specific spaces, which satisfy additional conditions (e.g., are independent of time etc.). Hence,  there cannot be a global solution branch in spaces where uniqueness is established while time-dependent external forces are just assumed to be in $L^2$. The argument holds only for dimension $n\geq 3$ (as the Leray projection term has not the same bound in lower dimension), and, in its present form, it makes no specific prediction about the type of singularity at the tip of the cone. 

\section{Comment to literature}

In \cite{T} it is shown that it seems unlikely that an argument based on the energy identity can succeed in proving uniqueness and global regular existence for the Navier Stokes equation. So it seems unlikely that an argument as cited in \cite{web} can succeed - at least if uniqueness is included in the statement. The argument presented here does not depend on the energy identity, and it shows that there are singular solution branches of the incompressible Euler equation.


\begin{thebibliography}{19}
\baselineskip=12pt


\bibitem{web}
	 {\em math.stackexchange.com/questions/634890 } 

\bibitem{T}
{\sc Tao, T.} {\em Finite time blowup for an averaged three-dimensional Navier Stokes equation}, Journal of the Am. Math. Soc., Jun. 2015.

\end{thebibliography}
\end{document}